\documentclass[11pt]{amsart}
\usepackage{bm}
\usepackage{url} 
\usepackage{verbatim}
\usepackage{geometry}
\usepackage[titletoc]{appendix}
\usepackage{amsmath,amsthm,amsfonts}
\usepackage{multirow}
\usepackage{nicefrac}
\usepackage{longtable}
\usepackage{color}
\usepackage{float}
\usepackage{graphicx, amssymb,graphics}
\usepackage{cite}
\usepackage{algorithm}
\usepackage{cancel}

\newtheorem{theorem}{Theorem}[section]

\newtheorem{lemma}[theorem]{Lemma}
\newtheorem{remark}[theorem]{Remark}

\theoremstyle{definition}

\numberwithin{equation}{section}

\makeatletter
\@addtoreset{equation}{section}
\makeatother

\newcommand{\R}{{\if mm {\rm I}\mkern -3mu{\rm R}\else \leavevmode\hbox{I}\kern -.17em\hbox{R} \fi}}

\newcommand{\cL}{\mbox{{${\mathcal{L}}$}}}
\newcommand{\cD}{\mbox{{${\mathcal{D}}$}}}

\newcommand{\cH}{\mbox{{${\mathcal{H}}$}}}
\title[]{Convergence of long-time stable variable-step arbitrary order ETD-MS scheme for gradient flows with Lipschitz nonlinearity}

\author[Wenbin Chen]{Wenbin Chen}

\address[Wenbin Chen]{Shanghai Key Laboratory for Contemporary Applied Mathematics, School of Mathematical Sciences, Fudan University, Shanghai, China 200433}
\email{{\tt wbchen@fudan.edu.cn}}

\author[Zhaohui Fu]{Zhaohui Fu}
\address[Zhaohui Fu]{Department of Mathematics, National University of Singapore, Singapore, 119076}
\email{\tt fuzhmath@nus.edu.sg}

\author[Shufen Wang]{Shufen Wang}
\address[Shufen Wang]{Shenwan Hongyuan Securities Co., Ltd., Shanghai, China 200031}
\email{\tt 17110180015@fudan.edu.cn}

\author[Xiaoming Wang]{Xiaoming Wang}
\address[Xiaoming Wang]{School of Mathematics Sciences, Eastern Institute of Technology, Ningbo, China 315200}
\email{\tt wxm.math@outlook.com}



\begin{document}

\begin{abstract}
    We analyze a variable-step extension of a family of arbitrarily high-order exponential time differencing multistep (ETD-MS) schemes recently developed by the authors. We prove that the schemes are unconditionally stable in the sense that a modified energy—representing a slight perturbation of the original energy—decreases monotonically over time, provided the nonlinearity is Lipschitz continuous in some appropriate sense. Moreover, we establish optimal-order convergence under mild conditions on the time-step size and local time-step ratio. Numerical experiments on the thin film epitaxial growth model without slope selection, employing a novel variable-step second-order scheme, validate the theoretical findings as well as its potential in developing highly efficient time-adaptive solution.
\end{abstract}
\maketitle

%%%%%%%%%%%%%%%%%%%%%%%%%%%%%%%%%%%%%%%%%%%%%%%%%%%%%%%%%%%%%%%%%%%%%%%%%%%%%%%%

\section{Introduction} \label{sec1}
The long-time behavior of many mathematical models is of great interest, particularly in the context of gradient flows arising in materials science, where physically important processes such as coarsening evolve over extended time scales. Capturing such dynamics efficiently poses significant computational challenges, thereby motivating the development of highly effective numerical methods.

Two key strategies are commonly employed to enhance computational efficiency:
	1.	High-order schemes, which enable the use of relatively large time steps while maintaining a prescribed accuracy, thus reducing the total number of steps required to reach a given final time.
	2.	Adaptive (variable) time-stepping, which allows time steps to be dynamically adjusted, taking larger steps when the system evolves slowly, and smaller steps when rapid changes occur.

A substantial body of literature exists on designing energy-stable schemes and conducting long-time simulations of various gradient flows in materials science and fluid dynamics. Notable approaches include convex splitting, truncation techniques, the scalar auxiliary variable (SAV) method, and the invariant energy quadratization (IEQ) method; see, for example, \cite{wang2010Unconditionally, shen2012second, elliott1993global, eyre1998unconditionally, shen2010numerical, shin2017unconditionally, shen2018scalar, shen2019new, gong2020arbitrarily} and references therein.

Exponential time differencing (ETD) is a particularly attractive class of time discretization techniques that achieves high-order accuracy through the exact treatment of the linear part and the application of Duhamel’s principle to the nonlinear part \cite{beylkin1998new, du2004stability, du2005analysis}. The resulting formulation introduces a nonlinear integral term, for which two primary approximation strategies are commonly adopted: Runge–Kutta (RK) methods \cite{cox2002exponential, hochbruck2005explicit, hochbruck2010exponential} and multistep (MS) methods \cite{hochbruck2010exponential, hochbruck2011exponential, hairer1993n}. These approximations are typically explicit in nature to preserve the computational efficiency of the ETD framework.

While it is relatively straightforward to formally construct variable step arbitrarily high-order ETD schemes via the RK or MS approaches, ensuring energy stability is far more delicate, primarily due to the explicit treatment of the nonlinear term. Among these approaches, the ETD-RK method has been extensively studied and applied, offering a flexible framework to achieve high-order accuracy, see \cite{FY2022ETDRK2JCP,CYC2024ETDRK32JCP,FSY2025ETDRKSCM,CYC2024,YCC2025 , quan2025highorder}. Recent advances also have addressed this issue by incorporating stabilization techniques into the ETD-MS framework, leading to energy-stable schemes of second order \cite{chen2020stabilized, ju2018energy}, third order \cite{cheng2019third, chen2020energy}, and arbitrary order \cite{chen2021ETDMS}, all under the assumption of a constant time step.

The primary objective of this paper is to extend the energy-stable ETD-MS framework to the variable time-step setting for a class of nonlinear gradient flows, building on our recent work in \cite{chen2021ETDMS}. We prove that the proposed variable-step ETD-MS schemes remain energy stable in the sense that a modified energy—a slight perturbation of the original energy—decreases monotonically over time. In addition, we establish that the schemes achieve optimal convergence rates under mild conditions on the time-step size and local time-step ratio.

Unless specified otherwise, we denote $C$ a generic constant which may depend on $\epsilon$, the exact solution $u$, the initial value $u_0$, and the final time $T$, but is independent of time step-size $\tau$. The standard Sobolev space notations follow from \cite{adams2003}.
	
The rest of the paper is organized as follows. The continuous problem is introduced in section 2. We introduce the variable-step arbitrary high-order ETD-MS based scheme and verify its energy stability in section 3. The convergence with optimal rate under a mild assumption on the step-size and the ratio of different step-sizes is established in section 4.
The applicability to the no-slope selection (NSS) thin-film epitaxial growth model is explained in section 5.  Numerical experiments on the no-slope selection thin film epitaxial growth model that are consistent with our theory are reported in section 6. Conclusion remarks are provided in section 7. 

\section{The continuous problem} \label{sec: continuous problem}
Let $X\hookrightarrow \cH\hookrightarrow X'$ be a standard triple of Hilbert spaces. And let $\cL$ be a positive definite linear operator on $\cH$ which maps $X$ to $X'$, the dual space of $X$.  Denoting $\cL^{1/2}$ the square root of $\cL$.
Let $E(u)$ be an energy functional on a Hilbert space $\cH$ (with the domain being a subspace of $\cH$) given by
$$E(u) = \frac{\epsilon}{2}\|\cL^{1/2} u\|_{\mathcal{H}}^2 + \mathcal{F}(u).$$
The associated gradient flow with mobility $M$ can then be formulated as follows
\begin{align} \label{eq:gradient flow}
\frac{\partial u}{\partial t}  = -M\frac{\delta E}{\delta u}.
\end{align}
where $\frac{\delta E}{\delta u}$ denotes the variational derivative of $E$.

Denoting the variational derivative of $\mathcal{F}$ as $-F$, i.e., $\frac{\delta \mathcal{F}}{\delta u} = -F $,
we deduce that the associated gradient flow takes the form of 
\begin{align} \label{eq:original eq}
\frac{\partial u(t)}{\partial t} + \epsilon \cL u(t)  = F(u(t)).
\end{align}
when $M=1$.

Taking the inner product of \eqref{eq:gradient flow} with $\frac{\delta E}{\delta u}$ in $\cH$ leads to the energy equality:
\begin{align} \label{eq:energy decay GF}
    \frac{dE(u)}{dt} = -M\left\|\frac{\delta E}{\delta u}\right\|_{\mathcal{H}}^2.
\end{align}

We will impose the following three assumptions on \eqref{eq:original eq}, with assumption $1-2$ for energy stability and the third one for convergence analysis:
\begin{enumerate}
    \item  The operator $\cL$ is self-adjoint nonnegative on $\mathcal{H}$. \\
    Thus we can define operators  $\cL^{\alpha/2}$ for any $\alpha\ge0$. The domain of operator $\cL^{\alpha/2}$ is denoted by $V^{\alpha} $ with the norm on $V^\alpha$ given by $\|v\|_{V^\alpha} = \|\cL^{\alpha/2}v\|_{\mathcal{H}}$. For $\alpha=1$ and $\alpha=0$,  it is abbreviated as $\cD(\cL^{\frac12}) = V$ and $\cD(\cL^0) = \cH$, respectively. 
    \item {The nonlinear term is Lipschitz continuous in the sense that: $\exists \beta\ge0, \gamma\ge0$, and $C_{L}>0$ with $\beta+\gamma\le 1$, such that }
    \begin{align} \label{eq:LIP}
    \left\| F(u) - F(v) \right\|_{V^{-\beta}} \le C_{L} \left\| u-v\right\|_{V^\gamma},
    \end{align}
    {where $V^{-\beta}$ is the dual space to $V^{\beta}$ induced by the inner product on $\cH$}.
%		\item Moreover, a specific Lipschitz continuity is required 
%		\begin{align} \label{eq:LIP2}
%		\left\| F(u) - F(v) \right\|_{V^{-1}} \le C_{L} \left\| u-v\right\|_{\mathcal{H}},
%		\end{align} 
%		with $C_{L}$ a positive constant.
    \item $F(u(\cdot,t))$ has the boundedness property in the sense that for any $T>0$, $\exists C=C(T)>0$, s.t.
    \begin{align} \label{eq:boundedness F}
    \left\| F(u) \right\|_{H^k(0,T;V^{-\beta})} \le C.
    \end{align}
\end{enumerate}

\begin{remark} \label{rem: condition for convergence}
    The restriction of $\beta+\gamma\le 1$ in the second assumption is for convergence only. 
\end{remark}
\begin{remark}	
    Likewise,  assumption 3 is used in the error estimate only. It can be verified for many systems such as the no-slope selection thin film epitaxial growth model investigated in the numerical experiment in this paper provided that the initial data is sufficiently smooth.
\end{remark}

\section{The temporal semi-discrete scheme and its long time stability}\label{sec: temporal discrete}	

The main purpose of this section is to introduce the variable step high-order ETD-MS based time discretization of the gradient flow \eqref{eq:original eq}, and prove the long time energy stability of the scheme. The scheme is a variable-step version of the arbitrary high-order ETD-MS scheme that was introduced in \cite{chen2021ETDMS}. The long time stability is in the sense that there is a modified energy that decreases in time. The modified energy is a small perturbation of the original one.

\subsection{The scheme}
	To solve \eqref{eq:original eq} numerically, we  propose the following ETD-MS based variable-step temporal semi-discrete numerical scheme: find $u^{n+1}(t)$ such that
	\begin{align} \label{eq:etdmsns}
	\frac{du^{n+1}(t)}{dt} + \epsilon \cL u^{n+1}(t) + A\tau^k \frac{d}{dt} \cL^{p(k)} u^{n+1}(t) = \sum_{i=0}^{k-1} \ell_i(t-t_n) F(u^{n-i}), \quad t\in[t_n,t_{n+1}],
	\end{align}
	with $u^{n+1}(t)=u^{n}(t)$ for $t\in[0,t_n]$. The temporal step-sizes are denoted by $\tau_{n}:=t_{n+1}-t_n$, and $\tau=\max\tau_j$.
	% and $\tau=\max\limits_{0\le j\le n}\tau_{j}$. 
	Moreover, we denote $u^{n-i}$  the numerical solution at time $t_{n-i}$ and $\ell_i(s)$  the shifted (to the negative range) Lagrange basis polynomial of degree $k$ with the form of
	\begin{align} \label{eq:Lagrange basis}
	\ell_i(s) = \prod_{0\le m\le k-1 \atop m\neq i} \frac{\sum_{j=1}^{m}\tau_{n-j}+s}{\sum_{j=1}^{m}\tau_{n-j}-\sum_{j=1}^{i}\tau_{n-j}} = \prod_{0\le m\le k-1 \atop m\neq i} \frac{\sum_{j=1}^{m}\frac{\tau_{n-j}}{\tau}+\frac{s}{\tau}}{\sum_{j=1}^{m}\frac{\tau_{n-j}}{\tau}-\sum_{j=1}^{i}\frac{\tau_{n-j}}{\tau}} = \sum_{j=0}^{k-1}\xi_{i,j}s^j,
	\end{align}
	where $\{ \xi_{i,j} \}_{j=0}^{k-1}$ are coefficients of the polynomial $\ell_{i}(s)$. Note that $\frac{\tau_j}{\tau}\le 1$ for any $j\le n$, it is not hard to see $\xi_{i,j} \sim \mathcal{O}(\tau^{-j})$. \footnote{We have suppressed the dependence of $\ell_i$ on $n$ for simplicity. The bound on $\ell_i$ is uniform in $n$ provided that the ratio between the maximum time-step and the minimum time-step in the neighboring $k$ steps is bounded by a constant independent of the step-sizes.}
	
	The scheme is the variable-step version of the scheme proposed recently by the authors \cite{chen2021ETDMS}.
	As in that work, we have employed a Dupont-Douglas type regularization term $A\tau^k \frac{d}{dt} \cL^{p(k)} u^{n+1}(t)$ in order to enhance stability. The parameters $A, p(k)$, associated with the strength of regularization, will be specified later.

%
%	Introducing an integrating factor $e^{\cK t}$ with $\cK=\epsilon\left(I+A\tau_n^k\cL^{p(k)}\right)^{-1}\cL$ and integrating the equation \eqref{eq:etdmsns} from $t_n$ to $t_{n+1}$ give the ETD-MS scheme
%	\begin{align} \label{eq:etdmsns 1}
%	u^{n+1} =& e^{-\cK \tau_n}u^n + \sum_{i=0}^{k-1}\int_{0}^{\tau} e^{-\cK (\tau_n -s)} \ell_i(s)  ds\, F(u^{n-i}).
%	\end{align}
%	Denote the integral $\int_{0}^{\tau} e^{-\cK (\tau_n -s)}  s^j $ by $\phi_j$, then it can be calculated by the recurrence formula
%	\begin{equation} \label{eq:phi recur}
%	\left\{
%	\begin{aligned} 
%	\phi_0 = & \cK^{-1}\left(I - e^{-\cK\tau_n}\right), \\
%	\phi_j = & \cK^{-1}\left(I -\frac{j}{\tau_n}\phi_{j-1} \right), \quad 1\le j \le k-1. 
%	\end{aligned} \right.
%	\end{equation}
%	Thus \eqref{eq:etdmsns 1} can be further rewritten as
%	\begin{align} \label{eq:etdmsns 2}
%	u^{n+1} = e^{-\cK \tau_n}u^n + \sum_{i=0}^{k-1}\sum_{j=0}^{k-1} \xi_{i,j}\phi_j\,F(u^{n-i}).
%	\end{align}

	\subsection{Energy stability} \label{sec31}
	In this subsection, we establish the energy stability for the scheme \eqref{eq:etdmsns}. The proof is similar to the stability proof  in the constant step-size case presented in \cite{chen2021ETDMS} adapted to the current variable-step scenario. 
	{The interpolation inequalities utilized below involving various exponents of $V$, i.e., the domain of various powers of $\cL$, follow from the spectral representation of the self-adjoint non-negative operator $\cL$ \cite{Lax2002book} and H\"older's inequality.}
	
	First we present a few interpolation estimates that will be used later.
	\begin{lemma} \label{lem:interpolation}
		Let $\beta,\gamma$ be two non-negative numbers, $q\in [0,1]$, and $p(k)$ is chosen so that $p(k)>\max\{\beta,\gamma\}$. Then for any $u\in V^{\beta}$ and $v\in V^{\gamma}$, and arbitrary positive constants $\hat{C}, \tilde{C}$, the following inequalities hold for different cases of $\beta$ and $\gamma$. 
		\begin{enumerate}
%{\color{blue}
		\item  If $ \min\{\beta,\gamma\} > 0$, we have
		\begin{align} \label{eq:bad term}
		\tau \left\| u\right\|_{V^\beta} \left\| v\right\|_{V^\gamma} \le C_1 \left\| u\right\|_{\mathcal{H}}^{2}  + C_2 \tau^{\frac{2qp(k)}{\beta}} \left\|u\right\|_{V^{p(k)}}^{2} + C_3 \left\| v\right\|_{\mathcal{H}}^{2} + C_4 \tau^{\frac{2(1-q)p(k)}{\gamma}} \left\| v\right\|_{V^{p(k)}}^{2},
		\end{align}
		where 
		\begin{equation} \label{eq:c1c2c3c4}
		\begin{aligned} 
		C_1(\hat{C}) &= \frac{1-\beta/p(k)}{2} \hat{C}^{1/(1-\beta/p(k))}, \quad C_2(\hat{C}) = \frac{\beta}{2p(k)}\hat{C}^{-p(k)/\beta}, \\
		C_3(\tilde{C}) &= \frac{1-\gamma/p(k)}{2} \tilde{C}^{1/(1-\gamma/p(k))}, \quad C_4(\tilde{C}) = \frac{\gamma}{2p(k)}\tilde{C}^{-p(k)/\gamma}.
		\end{aligned}
		\end{equation}
	\item 	If  $\min\{\beta, \gamma\}=0$ but $\max\{\beta, \gamma\}>0$, we have
		\begin{align} \label{eq:bad term 1}
			\tau \left\| u\right\|_{V^\beta} \left\| v\right\|_{V^\gamma} \le & \frac12 \left\| u\right\|_{\mathcal{H}}^{2}  + C_3 \left\| v\right\|_{\mathcal{H}}^{2} + C_4 \tau^{\frac{2p(k)}{\gamma}} \left\| v\right\|_{V^{p(k)}}^{2}, \quad\mbox{if}\ \beta=0,\ ~\gamma>0\\
			\tau \left\| u\right\|_{V^\beta} \left\| v\right\|_{V^\gamma} \le & C_1 \left\| u\right\|_{\mathcal{H}}^{2}  + C_2 \tau^{\frac{2p(k)}{\beta}} \left\|u\right\|_{V^{p(k)}}^{2} + \frac12 \left\| v\right\|_{\mathcal{H}}^{2} , \quad\mbox{if}\  \beta>0,\ ~\gamma=0.
		\end{align}
	\item 	If  $\max\{\beta, \gamma\}=0$, we have
		\begin{align} \label{eq:bad term 2}
		\tau \left\| u\right\|_{V^\beta} \left\| v\right\|_{V^\gamma} \le & \frac{\tau}{2} \left\| u\right\|_{\mathcal{H}}^{2}  + \frac{\tau}{2} \left\| v\right\|_{\mathcal{H}}^{2} .
		\end{align}  
		%}
	\end{enumerate}
	\end{lemma}
	\begin{proof}	
		For $\beta,\gamma>0$, 	we employ interpolation inequality to control $\|\cdot\|_{V^\beta}$, $\|\cdot\|_{V^\gamma}$ by $\|\cdot\|_{H}$ and $\|\cdot\|_{V^{p(k)}}$, and we have:
		\begin{align} \label{eq:interpolation}
		\tau \left\| u\right\|_{V^\beta} \left\| v\right\|_{V^\gamma} \le \tau \left\| u\right\|_{\mathcal{H}}^{1-\beta/p(k)}  \left\| u\right\|_{V^{p(k)}}^{\beta/p(k)} \cdot \left\| v\right\|_{\mathcal{H}}^{1-\gamma/p(k)}  \left\| v\right\|_{V^{p(k)}}^{\gamma/p(k)} .  
		\end{align}
		Denoting the right hand side (RHS) of \eqref{eq:interpolation} by $I_1$, and invoking  Young's inequality twice, we deduce
		\begin{align} \label{eq:bad term2}
		I_1 \le &  \frac12 \tau^{2q} \left\| u\right\|_{\mathcal{H}}^{2-2\beta/p(k)}  \left\| u\right\|_{V^{p(k)}}^{2\beta/p(k)}  + \frac12 \tau^{2-2q} \left\| v\right\|_{\mathcal{H}}^{2-2\gamma/p(k)}  \left\| v\right\|_{V^{p(k)}}^{2\gamma/p(k)} \nonumber \\
		\le & \frac{1-\beta/p(k)}{2}\left(\hat{C} \left\| u\right\|_{\mathcal{H}}^{2-2\beta/p(k)}  \right)^{1/(1-\beta/p(k))}  + \frac{\beta}{2p(k)}\left(  \frac{\tau^{2q}}{\hat{C}}\left\| u\right\|_{V^{p(k)}}^{2\beta/p(k)} \right)^{p(k)/\beta} \nonumber \\
		& + \frac{1-\gamma/p(k)}{2}\left( \tilde{C}\left\| v\right\|_{\mathcal{H}}^{2-2\gamma/p(k)}  \right)^{1/(1-\gamma/p(k))}  + \frac{\gamma}{2p(k)}\left( \frac{\tau^{2-2q}}{\tilde{C}} \left\| v\right\|_{V^{p(k)}}^{2\gamma/p(k)} \right)^{p(k)/\gamma} \nonumber \\
		= & C_1 \left\| u\right\|_{\mathcal{H}}^{2}  + C_2 \tau^{\frac{2qp(k)}{\beta}} \left\|u\right\|_{V^{p(k)}}^{2} + C_3 \left\| v\right\|_{\mathcal{H}}^{2} + C_4 \tau^{\frac{2(1-q)p(k)}{\gamma}} \left\| v\right\|_{V^{p(k)}}^{2}.
		\end{align}
	{ The proof  for the case of either or both $\beta$ and $\gamma$ are zero is  similar. This completes the proof of Lemma~\ref{lem:interpolation}. }
	\end{proof}
	
	\begin{remark}
		Note that  the estimates \eqref{eq:bad term 1}--\eqref{eq:bad term 2} are limit cases of \eqref{eq:bad term} with the constants $\hat{C}$ and $\tilde{C}$ chosen appropriately.  
\end{remark}
%	\begin{lemma} \label{lem:interpolation 2}
%		For any $q\in(0,1)$, $u\in V^{\beta}$ and $v\in V^{\gamma}$ with $ \max\{\beta,\gamma\} < 1$, and any constants $\hat{C}, \tilde{C}$, the following inequality holds
%		\begin{align} \label{eq:bad term 2}
%		\left\| u\right\|_{V^\beta} \left\| v\right\|_{V^\gamma} \le C_5 \left\| u\right\|_{\mathcal{H}}^{2}  + C_6  \left\|u\right\|_{V}^{2} + C_7 \left\| v\right\|_{\mathcal{H}}^{2} + C_8 \left\| v\right\|_{V}^{2},
%		\end{align}
%		where 
%		\begin{equation} \label{eq:c1c2c3c4 2}
%		\begin{aligned} 
%		C_5(\hat{C}) &= \frac{1-\beta}{2} \hat{C}^{1/(1-\beta)}, \quad C_6(\hat{C}) = \frac{\beta}{2}\hat{C}^{-1/\beta}, \\
%		C_7(\tilde{C}) &= \frac{1-\gamma}{2} \tilde{C}^{1/(1-\gamma)}, \quad C_8(\tilde{C}) = \frac{\gamma}{2}\tilde{C}^{-1/\gamma}.
%		\end{aligned}
%		\end{equation}
%		In particular, \eqref{eq:bad term 2} is also available for $\beta=1$ or $\gamma=1$ as long as we take $\hat{C}$ or $\tilde{C}$ to be zero.
%		%		{\color{red} and $\hat{C},~\tilde{C}$ are generic constants that will be determined in specific problem.}
%	\end{lemma}
%\begin{proof}
%	For $\max\{\beta,\gamma\}< 1$, the proof is similar to that of Lemma \ref{lem:interpolation}, and is omitted here. 
%\end{proof}

	Now we prove the energy stability. For simplicity, we denote $\|\cdot\|_{L^2(t_i,t_j;V^{\alpha})}$ by $\|\cdot\|_{L^2(I_{i}^{j-i};V^{\alpha})}$ hereafter.
	\begin{lemma} \label{lem:energy estimate1}
		For scheme \eqref{eq:etdmsns}, one has the following energy estimate.
		\begin{align} \label{eq:energy estimate 1}
		&E(u^{n+1}) - E(u^{n})+	\left\| 	\frac{du^{n+1}(t)}{dt} \right\|_{L^2(I_{n}^1;\mathcal{H})}^2 + A\tau^k \left\| 	\frac{du^{n+1}(t)}{dt}\right\|_{L^2(I_{n}^1;V^{p(k)})}^2  \nonumber \\
		\le & \sum_{j=0}^{k-1}C_L  \tau_{n-j}^{\frac12}\left\| 1- \sum_{i=-1}^{j-1}\ell_i(t-t_n) \right\|_{L^2(I_{n}^1)} \left\| \frac{du^{n-j+1}(t)}{dt} \right\|_{L^2(I_{n-j}^1;V^{\gamma})} \left\| \frac{du^{n+1}(t)}{dt} \right\|_{L^2(I_{n}^1;V^{\beta})},
		\end{align}
		where the convention $\ell_{-1}(t-t_n)=0$ has been used.
	\end{lemma}
	\begin{proof}
		To establish the desired energy estimates, we  subtract $F(u^{n+1}(t))$ from both sides of \eqref{eq:etdmsns} and take the inner product of the result with $	\frac{du^{n+1}(t)}{dt}$ on $\cH$, which gives
		\begin{align} \label{eq:energy stability1}
		&	\left\| 	\frac{du^{n+1}(t)}{dt} \right\|_{\mathcal{H}}^2 + A\tau^k \left\| 	\frac{du^{n+1}(t)}{dt}\right\|_{V^{p(k)}}^2 + \frac{d}{dt} E(u^{n+1}(t)) \nonumber \\
		=& \left( \sum_{i=0}^{k-1} \ell_i(t-t_n) F(u^{n-i}) - F(u^{n+1}(t)), \frac{du^{n+1}(t)}{dt}\right)_{\mathcal{H}}.
		\end{align}
		Integrating from $t_n$ to $t_{n+1}$ gives
		\begin{align} \label{eq:energy stability2}
		&	\left\| 	\frac{du^{n+1}(t)}{dt} \right\|_{L^2(I_{n}^1;\mathcal{H})}^2 + A\tau^k \left\| 	\frac{du^{n+1}(t)}{dt}\right\|_{L^2(I_{n}^1;V^{p(k)})}^2 + E(u^{n+1}) - E(u^{n}) \nonumber \\
		=& \int_{t_n}^{t_{n+1}}\left( \sum_{i=0}^{k-1} \ell_i(t-t_n) F(u^{n-i}) - F(u^{n+1}(t)), \frac{du^{n+1}(t)}{dt}\right)_{\mathcal{H}} dt.
		\end{align}
		
		Note that the sum of the Lagrange basis functions equals one, i.e., $\sum_{i=0}^{k-1} \ell_i(t-t_n) \equiv 1$. Thus terms within the integral in RHS of \eqref{eq:energy stability1} (denoted by NLT) can be rewritten as
		\begin{align} \label{eq:NLT}
		\text{NLT} = &  \left( \sum_{i=0}^{k-1} \ell_i(t-t_n)  \left( F(u^{n-i}) - F(u^{n+1}(t)) \right), \frac{du^{n+1}(t)}{dt}\right)_\mathcal{H} \nonumber \\	
		= &  \left( \sum_{i=0}^{k-1} \ell_i(t-t_n)  \left( F(u^{n-i}) - F(u^{n-i+1}) + \cdots+ F(u^{n})  -F(u^{n+1}(t)) \right), \frac{du^{n+1}(t)}{dt}\right)_\mathcal{H} \nonumber \\
		= & \left( \sum_{i=0}^{k-1} \ell_i(t-t_n)  \left( F(u^{n}) - F(u^{n+1}(t)) \right), \frac{du^{n+1}(t)}{dt}\right)_\mathcal{H} \nonumber \\
		& + \sum_{j=1}^{k-1}\left(  \sum_{i=j}^{k-1}\ell_{i}(t-t_n)  \left( F(u^{n-j}) - F(u^{n-j+1}) \right), \frac{du^{n+1}(t)}{dt}\right)_\mathcal{H} \nonumber \\
		= & \left(   F(u^{n}) - F(u^{n+1}(t)), \frac{du^{n+1}(t)}{dt}\right)_\mathcal{H} \nonumber \\ 
		& +\sum_{j=1}^{k-1} \left( \left(1- \sum_{i=0}^{j-1}\ell_i(t-t_n) \right)\left( F(u^{n-j}) - F(u^{n-j+1}) \right), \frac{du^{n+1}(t)}{dt}\right)_\mathcal{H}.
		\end{align}
		
		Utilizing the Cauchy-Schwarz inequality and the Lipschitz continuity assumption \eqref{eq:LIP}, we deduce
		\begin{align}  \label{eq:NLT1}
		&\int_{t_n}^{t_{n+1}} \left(   F(u^{n}) - F(u^{n+1}(t)), \frac{du^{n+1}(t)}{dt}\right)_\mathcal{H} dt \nonumber \\
		\le & C_{L} \int_{t_n}^{t_{n+1}} \left\| u^n -u^{n+1}(t) \right\|_{V^{\gamma}} \left\| \frac{du^{n+1}(t)}{dt} \right\|_{V^{\beta}} dt \nonumber \\
		\le & C_{L} \int_{t_n}^{t_{n+1}} \tau_n^{\frac12}\left\| \frac{du^{n+1}(t)}{dt} \right\|_{L^2(I_{n}^1;V^{\gamma})} \left\| \frac{du^{n+1}(t)}{dt} \right\|_{V^{\beta}} dt \nonumber \\
		\le & C_{L}  \tau_n \left\| \frac{du^{n+1}(t)}{dt} \right\|_{L^2(I_{n}^1;V^{\gamma})} \left\| \frac{du^{n+1}(t)}{dt} \right\|_{L^2(I_{n}^1;V^{\beta})} \nonumber \\
		:= & \mathrm{NLT_0},
		\end{align}
		where the last inequality follows from  H\"{o}lder's inequality.
		
		Similarly, the remaining terms in RHS of \eqref{eq:energy stability2} can be estimated as
		\begin{align} \label{eq:NLT2}
		&	\int_{t_n}^{t_{n+1}} \left(   \left(1- \sum_{i=0}^{j-1}\ell_i(t-t_n) \right) \left( F(u^{n-j}) - F(u^{n-j+1}) \right), \frac{du^{n+1}(t)}{dt}\right)_\mathcal{H} dt \nonumber \\
		\le & C_{L}  \tau_{n-j}^{\frac12} \left\| \frac{du^{n-j+1}(t)}{dt} \right\|_{L^2(I_{n-j}^1;V^{\gamma})} \left\| \frac{du^{n+1}(t)}{dt} \right\|_{L^2(I_{n}^1;V^{\beta})}  \left\| 1- \sum_{i=0}^{j-1}\ell_i(t-t_n) \right\|_{L^2(I_{n}^1)} \nonumber \\
		:= &\mathrm{NLT_j}, \quad 1\le j \le k-1.
		\end{align}
		This completes the proof.
	\end{proof}
	
	Next we give an upper bound on $ \left\| 1- \sum_{i=0}^{j-1}\ell_i(t-t_n) \right\|_{L^2(I_{n}^1)}$ and further provide the energy stability for our scheme \eqref{eq:etdmsns}. 
	Recall that in \eqref{eq:Lagrange basis}, the Lagrange basis $\ell_i(s)$ is expressed as the polynomial of $s$ with coefficients $\xi_{i,j}$. According to the properties of $\xi_{i,j}$, it's easy to see  
	\begin{align} \label{eq:ell L2 integral}
	\left\| 1- \sum_{i=0}^{j-1}\ell_i(t-t_n) \right\|_{L^2(I_{n}^1)} = 	\left\| 1- \sum_{i=0}^{j-1}\sum_{r=0}^{k-1} \xi_{i,r} (t-t_n)^r \right\|_{L^2(I_{n}^1)} \le C^{*}_j\tau^{1/2}, \quad 1\le j\le k-1, 
	\end{align}
	%{\color{red} 
	where the constants $C^{*}_j$ are independent of time step-size $\tau_{i}~(i\le n)$, $\tau$ and current time $t$. 
	
	For convenience,  we follow the convention of $C_0^*=1$ hereafter, and we introduce another sequence of positive constants $  \bm{\overline{C}}_j$ through
	\begin{align} \label{eq:cjbar}
	\bm{\overline{C}}_j := \sum_{r=0}^{k-1-j}C_{k-1-r}^*= \sum_{l=j}^{k-1}C_{l}^*.
	\end{align}
	It follows that
	\begin{equation}\label{eq:cjbar-rel}
	\bm{\overline{C}}_j = \bm{\overline{C}}_{j+1} +C_j^*, \quad \bm{\overline{C}}_{k-1} =C^*_{k-1}.
	\end{equation}
	We now introduce the following {\bf modified energy}
	\begin{align} \label{eq:modifed energy}
	\tilde{E}(u^n) =& E(u^n) + 
	 C_{L}C_3  \sum_{j=1}^{k-1} \bm{\overline{C}}_j \left\| \frac{du^{n-j+1}(t)}{dt}\right\|_{L^2(I_{n-j}^1;\mathcal{H})}^{2} %\nonumber\\& 
	 + C_{L} C_4 \sum_{j=1}^{k-1} \bm{\overline{C}}_j \tau^{k} \left\| \frac{du^{n-j+1}(t)}{dt}\right\|_{L^2(I_{n-j}^1;V^{p(k)})}^{2},
	\end{align}
	where both $C_3, ~C_4$  depend on $\hat{C}$, $\tilde{C}$ as specified  in \eqref{eq:c1c2c3c4}.  
\begin{remark}
%\color{blue}{
It is easy to see that the second term on the right hand side is of the order of $\tau$ while the third term is of the order of $\tau^2$. Therefore, the modified energy is a small perturbation of the original energy $E$ when the maximum time-step $\tau$ is small.
%}
\end{remark}
 	
	Thanks to \eqref{eq:c1c2c3c4},  $C_1, C_3$ can be made as small as we need so long as we set $\hat{C}, \tilde{C}$ small enough. Therefore, for 
	proper  constants $\hat{C}$, $\tilde{C}$,  and $A$, to be specified below,  the following inequalities hold
	\begin{align} 
	1 - C_{L} \left(C_3+C_1\bm{\overline{C}}_0 \right)  &\ge C_{L}C_3\bm{\overline{C}}_1  , \label{eq:constant choose1} \\
	 A - C_{L}\left(C_4+C_2\bm{\overline{C}}_0 \right) &\ge C_{L}C_4\bm{\overline{C}}_1. \label{eq:constant choose2}
	\end{align}
	Note that $C_0^{*}=1$, and hence $\bm{\overline{C}}_0=\bm{\overline{C}}_1+1$ according to \eqref{eq:cjbar-rel}. Therefore \eqref{eq:constant choose1}--\eqref{eq:constant choose2} are equivalent to
	\begin{align*} 
	 1 & \ge C_{L} \left(C_3+C_1 \right)\bm{\overline{C}}_0  , \\
	A &\ge C_{L}\left(C_4+C_2 \right)\bm{\overline{C}}_0 .
	\end{align*}
	This can be accomplished if we set  $\hat{C}$ and $\tilde{C}$ small enough so that
	\begin{align} \label{eq:c choose}
		(1-\beta/p(k))\hat{C}^{1/(1-\beta/p(k))} + (1-\gamma/p(k))\tilde{C}^{1/(1-\gamma/p(k))} \le 2/(C_{L}\bm{\overline{C}}_0),
	\end{align} 
	and then let
	\begin{align}\label{eq:A}
		A = &  C_{L}\left(\frac{\beta}{2p(k)}\hat{C}^{-p(k)/\beta}+\frac{\gamma}{2p(k)}\tilde{C}^{-p(k)/\gamma} \right)\bm{\overline{C}}_0.
	\end{align}
	We have \eqref{eq:constant choose1}--\eqref{eq:constant choose2} with the choice of $\hat{C}$, $\tilde{C}$, and $A$.
	
	We are now ready to prove the main result of the energy stability. 

	\begin{theorem} \label{thm:energy stability}
		Let  $\hat{C}$ and $\tilde{C}$be chosen so that  \eqref{eq:c choose} is satisfied, and the stabilized coefficient $A$ be specified in \eqref{eq:A}. Then the numerical scheme \eqref{eq:etdmsns} is energy stable in the sense that
		\begin{align} \label{eq:energy stable}
		\tilde{E}(u^{n+1}) \le \tilde{E}(u^n), \quad \forall n\ge k.
		\end{align}
	\end{theorem}
	\begin{proof}
		By \eqref{eq:ell L2 integral}, the estimate \eqref{eq:NLT2} for $\mathrm{NLT_{j}}, ~1\le j\le k-1$ can be simplified to
		\begin{align} \label{eq:NLT3}
		\mathrm{NLT_{j}} =  C_{L} C_j^* \tau_n^{\frac12}\tau_{n-j}^{\frac12} \left\| \frac{du^{n-j+1}(t)}{dt} \right\|_{L^2(I_{n-j}^1;V^{\gamma})} \left\| \frac{du^{n+1}(t)}{dt} \right\|_{L^2(I_{n}^1;V^{\beta})} .
		\end{align}
		Applying Lemma \ref{lem:interpolation} to \eqref{eq:NLT1} and \eqref{eq:NLT3}, these nonlinear terms can be bounded further: 
		\begin{align} 
		\mathrm{NLT_0} \le & C_{L} \left[ C_1 \left\| \frac{du^{n+1}(t)}{dt}\right\|_{L^2(I_{n}^1;\mathcal{H})}^{2}  + C_2 \tau_n^{\frac{2qp(k)}{\beta}} \left\|\frac{du^{n+1}(t)}{dt}\right\|_{L^2(I_{n}^1;V^{p(k)})}^{2} \right. \nonumber \\
		& \left. + C_3 \left\| \frac{du^{n+1}(t)}{dt}\right\|_{L^2(I_{n}^1;\mathcal{H})}^{2} + C_4 \tau_n^{\frac{2(1-q)p(k)}{\gamma}} \left\| \frac{du^{n+1}(t)}{dt}\right\|_{L^2(I_{n}^1;V^{p(k)})}^{2} \right],  \label{eq:interp1}\\
		\mathrm{NLT_{j}} \le & C_{L} C_j^*\left[ C_1\left\| \frac{du^{n+1}(t)}{dt}\right\|_{L^2(I_{n}^1;\mathcal{H})}^{2}  + C_2 (\tau_n\tau_{n-j})^{\frac{qp(k)}{\beta}} \left\|\frac{du^{n+1}(t)}{dt}\right\|_{L^2(I_{n}^1;V^{p(k)})}^{2} \right. \nonumber \\
		& \left. + C_3\left\| \frac{du^{n-j+1}(t)}{dt}\right\|_{L^2(I_{n-j}^1;\mathcal{H})}^{2} + C_4(\tau_n\tau_{n-j})^{\frac{(1-q)p(k)}{\gamma}} \left\| \frac{du^{n-j+1}(t)}{dt}\right\|_{L^2(I_{n-j}^1;V^{p(k)})}^{2} \right]. \label{eq:interp2}
		\end{align}	
	%	{\color{red} 
		We now pick $q$ and $p(k)$ in the following manner 
		\begin{equation}  \label{eq:p(k)2}
			\left\{ 
			\begin{aligned}
			q = \frac{1}{1+\gamma/\beta}, \quad	p(k) = \frac{(\beta+\gamma)k}{2}, & \quad \mbox{if}\ \beta>0, \gamma > 0, \\
			q =0, \quad	p(k) = \frac{(\beta+\gamma)k}{2}, & \quad \mbox{if}\ \beta=0, \gamma > 0, \\
			q =1, \quad	p(k) = \frac{(\beta+\gamma)k}{2}, & \quad \mbox{if}\ \beta>0, \gamma = 0. \\
			\end{aligned}
			\right.
		\end{equation} 
	We then have 
		\begin{equation} \label{eq:p(k)}
			\left\{ 
			\begin{aligned} 
			\frac{2qp(k)}{\beta} = k, \quad \frac{2(1-q)p(k)}{\gamma} = k, \quad & \quad \mbox{if}\ \beta>0, \gamma > 0, \\
			\  \quad \frac{2(1-q)p(k)}{\gamma} = k, \quad & \quad \mbox{if}\ \beta=0, \gamma > 0 ,\\
			\frac{2qp(k)}{\beta} = k, \quad \  \quad & \quad \mbox{if}\ \beta>0, \gamma = 0. \\
			\end{aligned}
			\right.
		\end{equation} 
		%}
		%{\color{red}
		Then estimates \eqref{eq:interp1}--\eqref{eq:interp2} give
		\begin{align} 
		\mathrm{NLT_0} \le & C_{L} \left[ C_1 \left\| \frac{du^{n+1}(t)}{dt}\right\|_{L^2(I_{n}^1;\mathcal{H})}^{2}  + C_2 \tau_n^{k} \left\|\frac{du^{n+1}(t)}{dt}\right\|_{L^2(I_{n}^1;V^{p(k)})}^{2} \right. \nonumber \\
		& \left. + C_3 \left\| \frac{du^{n+1}(t)}{dt}\right\|_{L^2(I_{n}^1;\mathcal{H})}^{2} + C_4 \tau_n^{k} \left\| \frac{du^{n+1}(t)}{dt}\right\|_{L^2(I_{n}^1;V^{p(k)})}^{2} \right],  \label{eq:interp3}\\
		\mathrm{NLT_{j}} \le & C_{L} C_j^*\left[ C_1\left\| \frac{du^{n+1}(t)}{dt}\right\|_{L^2(I_{n}^1;\mathcal{H})}^{2}  + C_2 (\tau_n\tau_{n-j})^{k/2} \left\|\frac{du^{n+1}(t)}{dt}\right\|_{L^2(I_{n}^1;V^{p(k)})}^{2} \right. \nonumber \\
		& \left. + C_3\left\| \frac{du^{n-j+1}(t)}{dt}\right\|_{L^2(I_{n-j}^1;\mathcal{H})}^{2} + C_4(\tau_n\tau_{n-j})^{k/2} \left\| \frac{du^{n-j+1}(t)}{dt}\right\|_{L^2(I_{n-j}^1;V^{p(k)})}^{2} \right]. \label{eq:interp4}
		\end{align}	
		%}	
		Simplifying the expression in \eqref{eq:interp1}--\eqref{eq:interp2} with  the convention of $C_0^*=1$ and combining \eqref{eq:interp3}--\eqref{eq:interp4}  with  \eqref{eq:energy stability2}, we have

		\begin{align} \label{eq:energy stability3}
		&	\left\| 	\frac{du^{n+1}(t)}{dt} \right\|_{L^2(I_{n}^1;\mathcal{H})}^2 + A\tau^k \left\| 	\frac{du^{n+1}(t)}{dt}\right\|_{L^2(I_{n}^1;V^{p(k)})}^2 + E(u^{n+1}) - E(u^{n}) \nonumber \\
		\le & C_{L} \left(C_3+C_1\sum_{j=0}^{k-1} C_j^* \right) \left\| \frac{du^{n+1}(t)}{dt}\right\|_{L^2(I_{n}^1;\mathcal{H})}^{2}  + C_{L}\left(C_4\tau_{n}^k+C_2\sum_{j=0}^{k-1} C_j^*(\tau_n\tau_{n-j})^{k/2} \right) \left\|\frac{du^{n+1}(t)}{dt}\right\|_{L^2(I_{n}^1;V^{p(k)})}^{2} \nonumber \\
		& + C_{L}C_3\sum_{j=1}^{k-1}  C_j^* \left\| \frac{du^{n-j+1}(t)}{dt}\right\|_{L^2(I_{n-j}^1;\mathcal{H})}^{2}  + C_{L}C_4 \sum_{j=1}^{k-1}  C_j^*(\tau_n\tau_{n-j})^{k/2}  \left\| \frac{du^{n-j+1}(t)}{dt}\right\|_{L^2(I_{n-j}^1;V^{p(k)})}^{2} \nonumber \\
		\le & C_{L} \left(C_3+C_1\sum_{j=0}^{k-1} C_j^* \right) \left\| \frac{du^{n+1}(t)}{dt}\right\|_{L^2(I_{n}^1;\mathcal{H})}^{2}  + C_{L}C_3\sum_{j=1}^{k-1}  C_j^* \left\| \frac{du^{n-j+1}(t)}{dt}\right\|_{L^2(I_{n-j}^1;\mathcal{H})}^{2}  \nonumber \\
		& + C_{L}\left(C_4\tau^k+C_2\sum_{j=0}^{k-1} C_j^*\tau^{k} \right)  \left\|\frac{du^{n+1}(t)}{dt}\right\|_{L^2(I_{n}^1;V^{p(k)})}^{2}   + C_{L}C_4 \sum_{j=1}^{k-1}  C_j^*\tau^{k}   \left\|\frac{du^{n-j+1}(t)}{dt}\right\|_{L^2(I_{n-j}^1;V^{p(k)})}^{2}  \nonumber \\
		= & C_{L} \left(C_3+C_1\bm{\overline{C}}_0 \right) \left\| \frac{du^{n+1}(t)}{dt}\right\|_{L^2(I_{n}^1;\mathcal{H})}^{2} + C_{L}C_3\sum_{j=1}^{k-1}  C_j^* \left\| \frac{du^{n-j+1}(t)}{dt}\right\|_{L^2(I_{n-j}^1;\mathcal{H})}^{2} \nonumber \\
		& + C_{L}\left(C_4+C_2\bm{\overline{C}}_0 \right) \tau^{k} \left\|\frac{du^{n+1}(t)}{dt}\right\|_{L^2(I_{n}^1;V^{p(k)})}^{2}   + C_{L}C_4 \sum_{j=1}^{k-1}  C_j^*\tau^{k}   \left\|\frac{du^{n-j+1}(t)}{dt}\right\|_{L^2(I_{n-j}^1;V^{p(k)})}^{2} ,
		\end{align}

		where the property of $\tau$ has been used.
		
		%{\color{red} 
		Adding  $C_{L}C_3 \sum_{j=1}^{k-2}\bm{\overline{C}}_{j+1} \left\| \frac{du^{n-j+1}(t)}{dt}\right\|_{L^2(I_{n-j}^1;\mathcal{H})}^{2}$ and  $C_{L}C_4 \sum_{j=1}^{k-2} \bm{\overline{C}}_{j+1}\tau^{k} \left\| \frac{du^{n-j+1}(t)}{dt}\right\|_{L^2(I_{n-j}^1;V^{p(k)})}^{2}$ to both sides of \eqref{eq:energy stability3}, and utilizing \eqref{eq:cjbar-rel}, we deduce
		\begin{align} \label{eq:energy stability5}
		&E(u^{n+1})+	\left( 1 - C_{L} \left(C_3 +C_1\bm{\overline{C}}_0  \right) \right)\left\| 	\frac{du^{n+1}(t)}{dt} \right\|_{L^2(I_{n}^1;\mathcal{H})}^2 \nonumber \\
		& + \left( A - C_{L}\left(C_4 +C_2\bm{\overline{C}}_0 \right) \right)\tau^k \left\| 	\frac{du^{n+1}(t)}{dt}\right\|_{L^2(I_{n}^1;V^{p(k)})}^2  \nonumber  \\
		& +  C_{L}C_3  \sum_{j=1}^{k-2} \bm{\overline{C}}_{j+1} \left\| \frac{du^{n-j+1}(t)}{dt}\right\|_{L^2(I_{n-j}^1;\mathcal{H})}^{2}  +  C_{L}C_4 \sum_{j=1}^{k-2} \bm{\overline{C}}_{j+1}\tau^{k} \left\| \frac{du^{n-j+1}(t)}{dt}\right\|_{L^2(I_{n-j}^1;V^{p(k)})}^{2} 
		\nonumber \\
		\le &  E(u^{n}) +  C_{L}C_3\sum_{j=1}^{k-1}  C_j^* \left\| \frac{du^{n-j+1}(t)}{dt}\right\|_{L^2(I_{n-j}^1;\mathcal{H})}^{2}  + C_{L}C_4 \sum_{j=1}^{k-1}  C_j^*\tau^{k}   \left\|\frac{du^{n-j+1}(t)}{dt}\right\|_{L^2(I_{n-j}^1;V^{p(k)})}^{2}   \nonumber \\
		&+  C_{L}C_3  \sum_{j=1}^{k-2} \bm{\overline{C}}_{j+1} \left\| \frac{du^{n-j+1}(t)}{dt}\right\|_{L^2(I_{n-j}^1;\mathcal{H})}^{2}  +  C_{L}C_4 \sum_{j=1}^{k-2} \bm{\overline{C}}_{j+1}\tau^{k} \left\| \frac{du^{n-j+1}(t)}{dt}\right\|_{L^2(I_{n-j}^1;V^{p(k)})}^{2}  \nonumber \\
		= &  E(u^{n}) +  C_{L}C_3 \bm{\overline{C}}_{1} \left\| \frac{du^{n}(t)}{dt}\right\|_{L^2(I_{n-1}^1;\mathcal{H})}^{2}  + C_{L}C_4 \bm{\overline{C}}_{1} \tau^{k} \left\| \frac{du^{n}(t)}{dt}\right\|_{L^2(I_{n-1}^1;V^{p(k)})}^{2}  \nonumber \\
		& +  C_{L}C_3  \sum_{j=2}^{k-1} \bm{\overline{C}}_j \left\| \frac{du^{n-j+1}(t)}{dt}\right\|_{L^2(I_{n-j}^1;\mathcal{H})}^{2}  + C_{L} C_4 \sum_{j=2}^{k-1} \bm{\overline{C}}_j \tau^{k} \left\| \frac{du^{n-j+1}(t)}{dt}\right\|_{L^2(I_{n-j}^1;V^{p(k)})}^{2} .
		\end{align}
		Then the modified energy decaying property \eqref{eq:energy stable}
		follows from \eqref{eq:constant choose1}--\eqref{eq:constant choose2} and \eqref{eq:energy stability5}.
		
	\end{proof}

%{\color{red}	

\begin{remark}
    {{For the case of $\beta=\gamma=0$, the nonlinear terms in \eqref{eq:NLT} can be bounded by  $C_0\tau_{n}\left\|\frac{du^{n+1}}{dt} \right\|^2_{L^2(I_n^1;\mathcal{H})} + \sum_{j=1}^{k-1}\tau_n\left\| \frac{du^{n-j+1}}{dt}\right\|^2_{L^2(I_{n-j}^1;\mathcal{H})}$ where we have used the boundedness of step ratios of neighboring $k$-steps. Therefore, we can define a modified energy of the form of 
    $\tilde{E}(u^n) = E(u^n) + 
	  \sum_{j=1}^{k-1}\tilde{C}_j \tau_{n}\left\| \frac{du^{n-j+1}(t)}{dt}\right\|_{L^2(I_{n-j}^1;\mathcal{H})}^{2} %\nonumber\\& 
	 $ with $\tilde{C}_{j+1}r_c+1\le \tilde{C}_{j}$ where $r_c$ is a bound on neighboring $k$ step sizes. %, 1\ge (\tilde{C}_1 r_c+ C_0)\tau_n, \tau_n\le \tilde{C}_{k-1}$.  
	 Such non-negative choices of $\tilde{C}_j$s are always possible. We can demonstrate the monotonic decreasing property of the modified energy without the regularization term in the energy equality (3.12) provided that the time-step is small ( $(C_0+\tilde{C}_1 r_c)\tau_n \le 1, \tau_n\le \tilde{C}_{k-1}$).  } }Thus no additional regularization term is required {for energy stability}, but a constant restriction for time step-size $\tau$ is needed.
\end{remark}

    %%%%%
	\section{Convergence analysis} \label{sec:conv}	
	In this section, we provide the optimal error estimate for the temporal discrete scheme  \eqref{eq:etdmsns} on any finite time interval $[0, T]$ assuming that the solution is sufficiently smooth so that assumption 3 is satisfied.
	
	Let $u(t)$ the exact solution of \eqref{eq:original eq} and define the error function $e(t)=u(t)-u^{n+1}(t), e^j=e(t_j)=u(t_j)-u^j$, then the error equation becomes
	\begin{align} \label{eq: error eq}
		& \frac{d e(t)}{dt} + \epsilon \cL e(t) + A\tau^k \frac{d}{dt}\cL^{p(k)} e(t) \nonumber \\
		= &  A\tau^k \frac{d}{dt}\cL^{p(k)} u(t) +  F(u(t)) -  \sum_{i=0}^{k-1} \ell_i(t-t_n) F(u^{n-i})  \nonumber \\
		:= & R_1 +  NL_1 + NL_2,
	\end{align}	
	where 
	\begin{align*}
		R_1 = &  A\tau^k \frac{d}{dt}\cL^{p(k)} u(t), \\
		NL_1 = & F(u(t)) -  \sum_{i=0}^{k-1} \ell_i(t-t_n)  F(u(t_{n-i})), \\
		NL_2 = & \sum_{i=0}^{k-1} \ell_i(t-t_n)  F(u(t_{n-i})) -  \sum_{i=0}^{k-1} \ell_i(t-t_n)  F(u^{n-i}).
	\end{align*}
	Assume
	\begin{equation}
	 \lambda=1-\beta. \label{eq: lambda}
	\end{equation}
	Taking the inner product of  \eqref{eq: error eq} with $\cL^\lambda e(t)$ in $\cH$ , we have for some generic constant $C$,
	\begin{align} \label{eq: error eq product}
	&  \frac12 \frac{d}{dt} \|e(t)\|^2_{V^\lambda} + \epsilon \|e(t)\|^2_{V^{\lambda+1}} + \frac12 A\tau^k \frac{d}{dt} \|e(t)\|^2_{V^{p(k)+\lambda}}  \nonumber \\
	=& \left(R_1+NL_1+NL_2,\cL^\lambda e(t)\right)_{\mathcal{H}}  \nonumber \\
	 \le &  C \|R_1\|^2_{V^{\lambda-1}} +  C \|NL_1\|^2_{V^{\lambda-1}} + \frac{\epsilon}{4} \|e(t)\|^2_{V^{\lambda+1}}  + \left( NL_2, \cL^\lambda e(t)\right)_{\mathcal{H}} .
	\end{align}
	It is easy to see 
	\begin{align} \label{eq:R1}
		\|R_1\|^2_{V^{\lambda-1}} \le C\tau^{2k}\left\| \frac{d u}{dt} \right\|_{V^{2p(k)+\lambda-1}}^2.
	\end{align}
	For $NL_1$, the properties of Lagrange interpolation give
	\begin{align} \label{eq: NL1}
		NL_1 = & \sum_{i=0}^{k-1} \ell_i(t-t_n) \left[  F(u(t)) - F(u(t_{n-i}))  \right] \nonumber \\
		= & \frac{1}{k!} \sum_{i=0}^{k-1} \ell_{i}(t-t_n) \int_{t}^{t_{n-i}} \left( t_{n-i} - t \right)^{k-1} \left( F(u(s))\right)_s^{(k)} ds,
	\end{align}
	where $\left( F(u(s))\right)_s^{(k)}$ represents the $k$-th derivative  with respect to $s$. Applying H\"{o}lder inequality to \eqref{eq: NL1} , we have
	\begin{align} \label{eq: NL1 norm}
	\| NL_1 \|^2_{V^{\lambda-1}} \le & C\tau^{2k-1}  \left\| F(u(s))\right\|^2_{H^k(I_{n-k+1}^{k};V^{\lambda-1})}.
	\end{align}
	Next, we turn to the other nonlinear term, $NL_2$. By \eqref{eq:LIP} one has
	\begin{align} \label{eq: NL2}
	  \left( NL_2,\cL^\lambda e(t)\right)_{\mathcal{H}} = & \sum_{i=0}^{k-1} \ell_{i}(t-t_n) \left( F(u(t_{n-i})) - F(u^{n-i}), \cL^\lambda e(t) \right)_{\mathcal{H}} \nonumber \\
	  \le & \sum_{i=0}^{k-1} |\ell_{i}(t-t_n)| \left\| F(u(t_{n-i})) - F(u^{n-i}) \right\|_{V^{-\beta}} \| e(t) \|_{V^{2\lambda+\beta}} \nonumber \\
	  \le & C_{L} \sum_{i=0}^{k-1}| \ell_{i}(t-t_n)| \left\| e^{n-i}  \right\|_{V^\gamma} \| e(t) \|_{V^{2\lambda+\beta}} \nonumber \\
	  \le & {\sum_{i=0}^{k-1} C\,C_L^2 \|e^{n-i}\|^2_{V^\gamma} + \frac{\epsilon}{4} \|e(t)\|^2_{V^{2\lambda+\beta}} }.%\nonumber \\
%	  = & \sum_{i=0}^{k-1}  \frac{1}{8\epsilon}\|e^{n-i}\|^2_{V^\gamma} +2C_{L}^2\epsilon \|e(t)\|^2_{V^{2\lambda+\beta}}.
	\end{align}	
	%{\color{blue}
	Combining \eqref{eq: error eq product},  \eqref{eq: NL1 norm} and \eqref{eq: NL2}, we obtain
	\begin{align} \label{eq: error eq product 2}
		 & \frac12 \frac{d}{dt} \|e(t)\|^2_{V^\lambda} + \frac{3\epsilon}{4} \|e(t)\|^2_{V^{\lambda+1}} + \frac12 A\tau^k \frac{d}{dt} \|e(t)\|^2_{V^{p(k)+\lambda}}  \nonumber \\
		\le & 
		 C\, C_L^2 \sum_{i=0}^{k-1} \|e^{n-i}\|^2_{V^\gamma} +\frac{\epsilon}{4} \|e(t)\|^2_{V^{2\lambda+\beta}} \nonumber \\
		&+ C\tau^{2k}\left\| \frac{d u}{dt} \right\|_{V^{2p(k)+\lambda-1}}^2  + C\tau^{2k-1} \left\| F(u) \right\|^2_{H^k(I_{n-k+1}^{k};V^{\lambda-1})}.
	\end{align}
	Note that $\beta = 1-\lambda$, then $2\lambda+\beta = \lambda + 1$. Denote $w(t)=\|e(t)\|^2_{V^\lambda} +  A\tau^k \|e(t)\|^2_{V^{p(k)+\lambda}}$, then \eqref{eq: error eq product 2} can be written as
	\begin{align} \label{eq: error eq product 3}
	 \frac{d}{dt} w(t)  + \epsilon \|e(t)\|^2_{V^{\lambda+1}}
	\le & C\,C_{L}^2 \sum_{i=0}^{k-1}  \left\| e^{n-i}  \right\|^2_{V^\gamma} 
	  + C\tau^{2k}\left\| \frac{du}{dt}\right\|_{V^{2p(k)-\beta}}^2  +
	C\tau^{2k-1} \left\| F(u) \right\|^2_{H^k(I_{n-k+1}^{k};V^{-\beta})}.
	\end{align}
	
	%{\color{blue}
	Since $\gamma \le 1 - \beta = \lambda$, we have $\|e^{n-i}\|^2_{V^\gamma} \le w^{n-i}$. Integrating  \eqref{eq: error eq product 3} from $t_n$ to $t_{n+1}$ , we deduce
	%}
	 	\begin{align} \label{eq: error eq product 5}
	w(t_{n+1}) - w(t_{n}) 
	 \le &C\,C_{L}^2 \sum_{i=0}^{k-1} w^{n-i} \tau_n dt  +  C\tau^{2k}\left\| \frac{du}{dt}\right\|_{L^2(I_n^1, V^{2p(k)-\beta})}^2 +
	 C\tau^{2k} \left\| F(u) \right\|^2_{H^k(I_{n-k+1}^{k};V^{-\beta})} .
	 \end{align}
	
	 Summing up for $n$ from $k-1$ to $m$, utilizing the assumptions on $\frac{du}{dt}$ and $F(u)$ with $t_{m+1}\le T$, we have
	  \begin{align} \label{eq: error eq product 8}
	 w^{m+1} - w^{k-1} %-  \sum_{n=k}^{m+1}w^{n}\tau_{n-1}
	 \le & C\,C_L^2\sum_{n=k-1}^m \sum_{i=0}^{k-1} w^{n-i}\tau_n  + C\tau^{2k}.
	 \end{align}
	 In order to apply the discrete Gronwall's inequality in summation form, we impose the following  time-step ratio on neighboring $k$ steps ({\bf local time-step ratio}). More specifically, we assume there exists $r_c>0$, such that
	\begin{eqnarray} 
	% 	 \tau_n &< & 1, \forall n, \label{eq: tau restriction}\\
	 	 \frac{\tau_n}{\tau_m} &\le& r_c, \quad \forall n,m, \ s.t.\ |n-m|<k. \label{eq: local step ratio restriction}
	 \end{eqnarray}
	  If we further assume that the initial errors are of the order of $k$, i.e.,  $w^{j} \le C\tau^{2k}, j=0.\cdots,k-1$, 
	 we deduce from \eqref{eq: error eq product 8}
	 \begin{equation}
	w^{m+1} 
	 \le  C\,C_L^2\sum_{n=1}^{m} w^{n}\tau_{n-1}  + C\tau^{2k}.
	 \end{equation}
	
	 We have, thanks to discrete Gronwall inequality
	 \begin{align} \label{eq: error end}
	 	w^m \le  C\tau^{2k}, \forall t_m \in [0,T].
	 \end{align}
	%}
    Therefore, we have proved the following optimal error estimates 
    % {\color{blue}
    \begin{theorem}
        For $\lambda=1-\beta$. Assume $F(u) \in H^k(0,T;V^{-\beta})$ and $u \in H^1(0,T;V^{2p(k)-\beta})$, $	\|e^{j}\|^2_{V^\lambda} +  A\tau^k \|e^{j}\|^2_{V^{p(k)+\lambda}} \le  C\tau^{2k}, j=0,1, \cdots, k-1$;  $\beta+\gamma \le 1$, $\tau < 1$, and the local  time-step ratio is bounded by a constant $r$ independent of the time step $\tau_n$, i.e., \eqref{eq: local step ratio restriction} is valid. Then  the scheme converges with the optimal rate of $k$ in the sense that 
        \begin{align} \label{eq: error estimate}
            \|e^n\|^2_{V^\lambda} +  A\tau^k \|e^n\|^2_{V^{p(k)+\lambda}} \le  C\tau^{2k}, \quad \forall n \quad s.t. \ t_{n}\le T.
        \end{align}
    \end{theorem} 
    %}

\section{Application to no-slope-selection (NSS) thin-film epitaxial growth model} \label{sec:application}

In this section, weshow that the abstract framework applies to the no-slope-selection (NSS) thin-film epitaxial growth equation:
\begin{align} \label{eq:nss}
\frac{\partial u}{\partial t}=-\epsilon \Delta^{2} u-\nabla \cdot\left(\frac{\nabla u}{1+|\nabla u|^{2}}\right)
\end{align}
with periodic boundary condition imposed. The energy functional is given by 
\begin{align} \label{eq:SS energy functional}
E(u)=\int_{\Omega} \left( \frac{\epsilon}{2}|\Delta u|^{2}-\frac{1}{2} \ln \left(1+|\nabla u|^{2}\right)  \right) \mathrm{d} \mathbf{x}.
\end{align}
The linear and nonlinear operators are $\cL = \Delta^2$ and $F(u) = -\nabla\cdot \left( \frac{\nabla u}{1+|\nabla u|^2}\right)$, respectively.
Therefore, the abstract functional spaces are specified to $\cH=\{f:f\in L^2~\text{with zero mean}\}$, $V^{\frac12}=\{f:f\in H^2_{per}~\text{with zero mean}\}$.

In the following, we will verify the assumptions proposed in Section \ref{sec: continuous problem} for the NSS equation one by one. First we note that the norm equivalence between $\|v\|_{H^2}$ and $\|\Delta v\|$ follows from the elliptic regularity and Poincaré's inequality. Secondly,  the Lipschitz continuity \eqref{eq:LIP} with $\beta=\gamma=\frac12$, $C_L=1$ has been verified in \cite{ju18energy}, i.e.,
\begin{align} \label{eq:LIP nss}
\left\| F(u) - F(v) \right\|_{V^{-\frac12}} \le&  \left\| u-v\right\|_{V^\frac12}.
\end{align}
By \eqref{eq:p(k)2}, it yields $p(k)=k/2$. Therefore, the modified energy stability \eqref{eq:energy stable} establishes with
\begin{align*}
	\tilde{E}_N(u^n) = &\left(-\frac12\ln\left(1+|\nabla u^n|^2\right),1\right)_{L^2} + \frac{\epsilon}{2}\|\Delta u^n\|_{L^2}^2 + C_L C_3  \sum_{j=1}^{k-1} \bm{\overline{C}}_j \left\| \frac{du^{n-j+1}(t)}{dt}\right\|_{L^2(I_{n-j}^1;L^2)}^{2} \nonumber\\
	& + C_LC_4 \sum_{j=1}^{k-1} \bm{\overline{C}}_j \tau^{k} \left\| \frac{d\Delta^{k/2}u^{n-j+1}(t)}{dt}\right\|_{L^2(I_{n-j}^1;L^2)}^{2}.
\end{align*}
Constructing proper initial step (the first $k-1$ steps) similar to that in \cite{chen2020energy} for solving \eqref{eq:nss},  then we can obtain an upper bound for $\tilde{E}(u^n)$, from which the $H^2$ bound of numerical solution to \eqref{eq:etdmsns} can be obtained (see \cite{chen2014linear}). Moreover, the $H^3$ bound for the numerical solution can be provided as long as the initial data is smooth enough (see \cite{chen2020energy,chen2020stabilized}). As for the boundedness \eqref{eq:boundedness F}, now $\lambda = 1 - \beta = \frac12$ and 
\begin{align} \label{eq:boundedness F proof}
	\left\| F(u) \right\|_{H^k(0,T;H^1)}^2 = & \left\|-\nabla \cdot \left( \frac{\nabla u}{1+|\nabla  u|^2}\right)\right\|_{H^k(0,T;H^1)}^2  \le  C\left\|  u \right\|_{W^{k,\infty}(0,T;H^3)}^2 .
\end{align}
thus this boundedness can be guaranteed by requiring exact solution $u$ to be smooth enough, say $u\in W^{k,\infty}(0,T;H^{m+2})$.

%%%%%
\section{Numerical results} \label{eq: numer}
In this section, we report numerical results when the  variable-step second-order ETD-MS scheme on the NSS equation. More specifically, the scheme takes the following form
\begin{align} \label{eq:second order nss}
& \frac{du^{n+1}(t)}{dt} + \epsilon\Delta^2 u^{n+1}(t) + A\tau^2 \frac{d}{dt} \Delta^2 u^{n+1}(t) =  -\nabla\cdot\left(\frac{\nabla u^n}{1+|\nabla u^n|^2}\right) \nonumber \\
&\quad + \frac{t - t_n}{\tau_n} \left(-\nabla\cdot\left(\frac{\nabla u^n}{1+|\nabla u^n|^2}\right) + \nabla\cdot\left(\frac{\nabla u^{n-1}}{1+|\nabla u^{n-1}|^2}\right) \right), \quad t\in[t_n,t_{n+1}].
\end{align}
Both the temporal convergence and long-time energy stability are validated. The two dimensional domain $\Omega=[0,2\pi]^2$ with periodic boundary condition is considered and the Fourier pseudo-spectral method is applied for  spatial discretization.

\subsection{Temporal convergence} \label{eq: time converge}
In this subsection, the second-order temporal convergence is tested at  the terminal time $T=1$.  The parameters are set as $L=4\pi$, $\epsilon=0.01$. 
The construction of variable time steps is completed via a uniform time grid  plus 10\% random perturbation; the specific  implementation is referred to \cite[Section 6.1, p.518]{Chen2019second}.
More specifically, the coarsest grid is obtained by setting a uniform partition with step-size $\Delta t_0=0.0025$, and then adding a $10\%$ random perturbation onto time grids to make  new variable-step time-step  series. The finer one is to double the number of time nodes with  the nodes  set as $t_k^{fine}=t_{(k+1)/2}^{coarse}$ for odd $k$, and $t_k^{fine}=\left( t_{k-1}^{fine} + t_{k+1}^{fine}\right)/2$ for even $k$. 

Simple calculation shows $C_0^*=1$, $C_1^*=1/\sqrt{3}$ and $\bm{\overline{C}}_0 =1+1/\sqrt{3}$. The constants in \eqref{eq:c1c2c3c4} are $C_1=\frac14\hat{C}^2, \ C_2=\frac14\hat{C}^{-2}, \ C_3=\frac14\tilde{C}^2, \ C_4=\frac14\tilde{C}^{-2}$. Then we can take $\hat{C}^2=\tilde{C}^2=\frac{6}{3+\sqrt{3}}$ and $A= \left(\frac{3+\sqrt{3}}{6}\right)^2 = \frac{2+\sqrt{3}}{6}$.

To compute the error, an artificial forcing term is added so that the exact solution is given by $u(t)=\cos(t)\sin(x)\cos(y)$.   Table~\ref{tab1: error} exhibits the $L^2$ error, in which the second-order accuracy is validated. And we can find that the errors on the uniform time mesh are smaller than those on the nonuniform mesh.

\begin{table}
    \centering
    \begin{tabular}{|c||c|c||c|c|}\hline
     $N_T$    &$L^2$ error  & convergence rate& $L^2$ error  & convergence rate \\\hline
      1 &5.258e-2&-& 7.873e-2 & -   \\\hline
     2 &1.658e-2&1.665 & 2.801e-2 & 1.4909   \\\hline
      4 &4.432e-3 &1.903& 7.807e-3   & 1.843 \\\hline
       8  &1.128e-3 &1.973& 2.011e-3    &1.957 \\\hline
       16&2.837e-4&1.992& 5.071e-4  &1.988 \\ \hline
 32   &7.103e-5&1.998& 1.271e-4 & 1.996\\ \hline
        {64} &1.777e-5&1.999&3.181e-5  & 1.999  \\ \hline
    \end{tabular}
    \caption{$L^2$ errors on  uniform initial time grid(the second column) and random initial time mesh(the fourth column), the time step $\Delta t=\Delta t_0/N_T$ for the uniform time grid.}
    \label{tab1: error} 
\end{table}

\subsection{Coarsening process} \label{subsec: energy}
\begin{figure}[ht]
	\centering
	 \noindent\makebox[\textwidth][c] {
		\begin{minipage}{0.3\textwidth}
	 		\includegraphics[width=\textwidth]{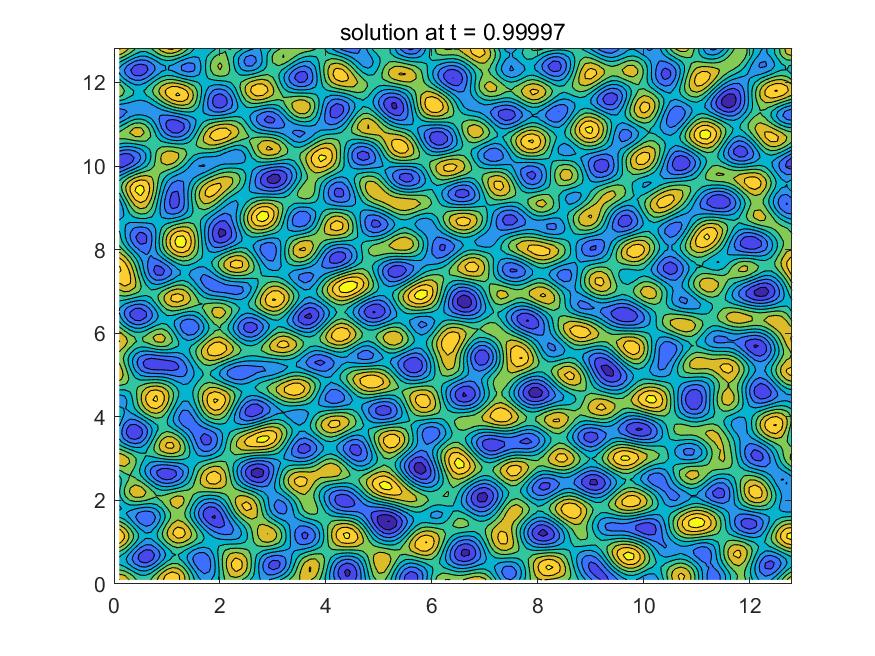}
		\end{minipage}
	 	\begin{minipage}{0.3\textwidth}
			\includegraphics[width=\textwidth]{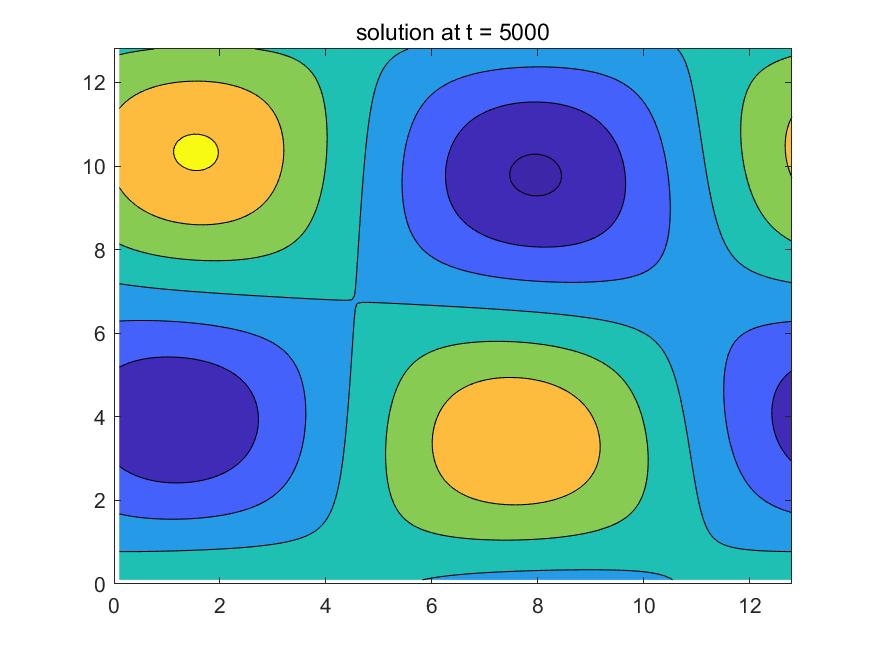}
	 	\end{minipage}
	 	\begin{minipage}{0.3\textwidth}
			\includegraphics[width=\textwidth]{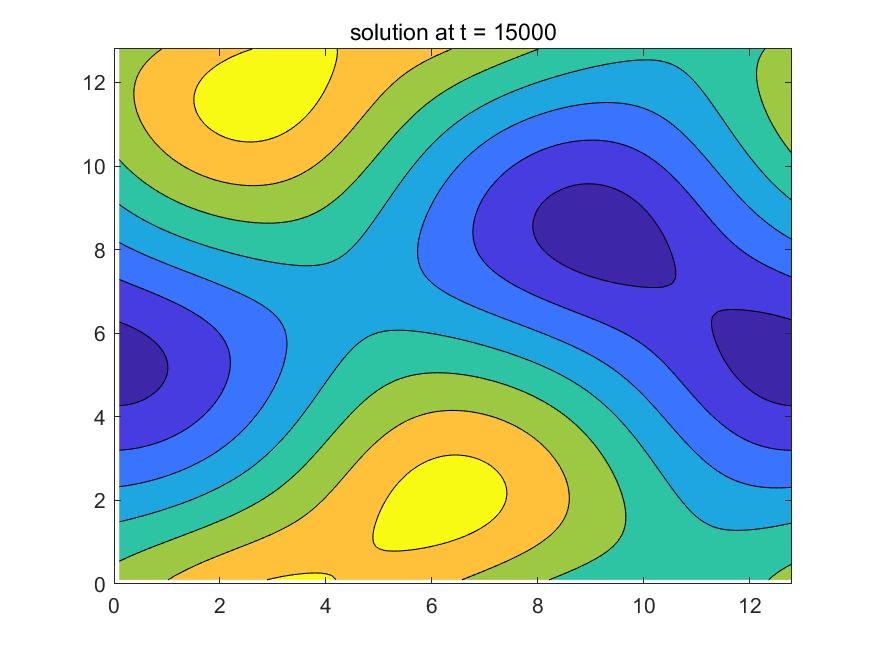}
	 	\end{minipage}
	 }
	 \noindent\makebox[\textwidth][c] {
		\begin{minipage}{0.3\textwidth}
	 		\includegraphics[width=\textwidth]{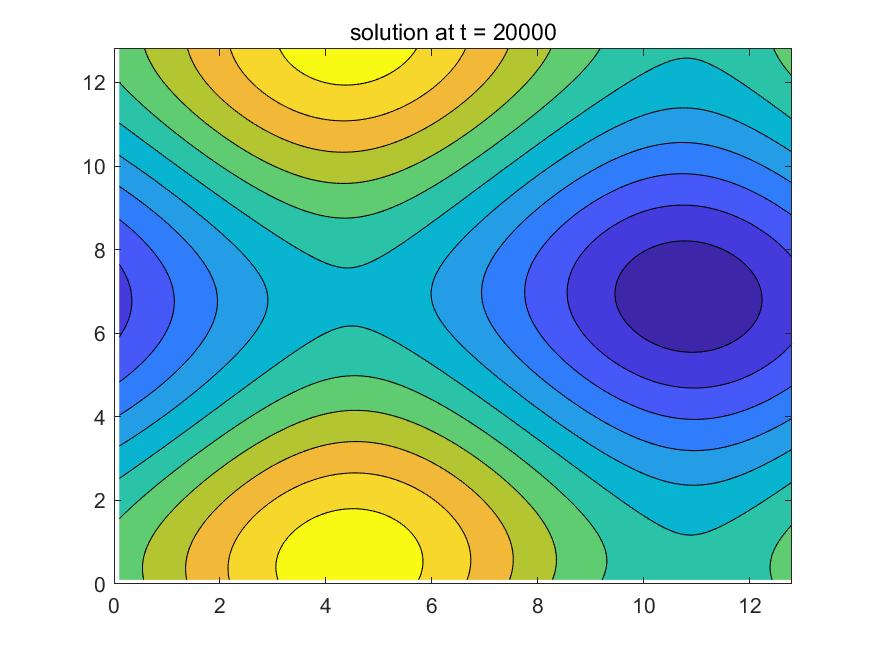}
		\end{minipage}
	 	\begin{minipage}{0.3\textwidth}
	 		\includegraphics[width=\textwidth]{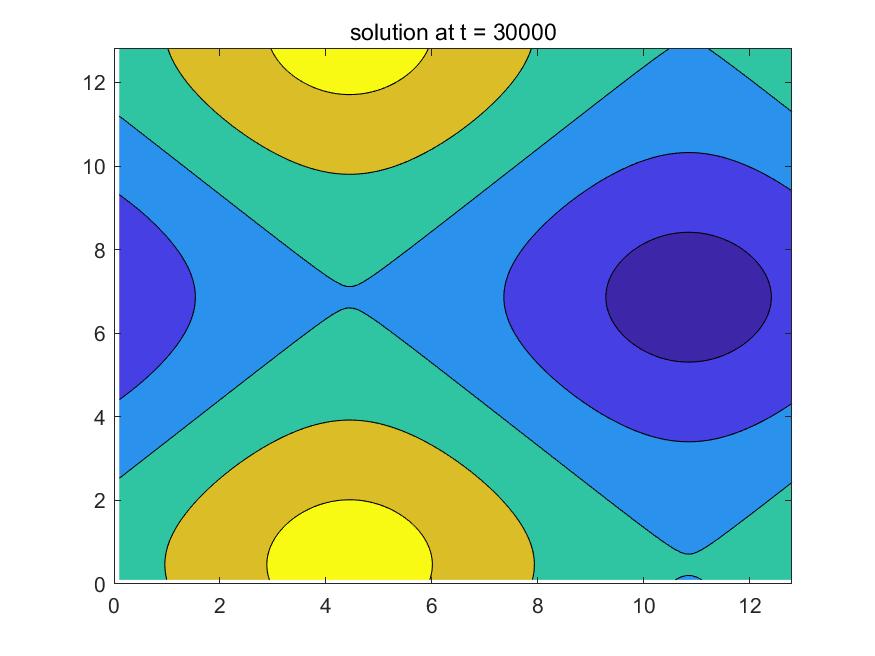}
	 	\end{minipage}
	 	\begin{minipage}{0.3\textwidth}
			\includegraphics[width=\textwidth]{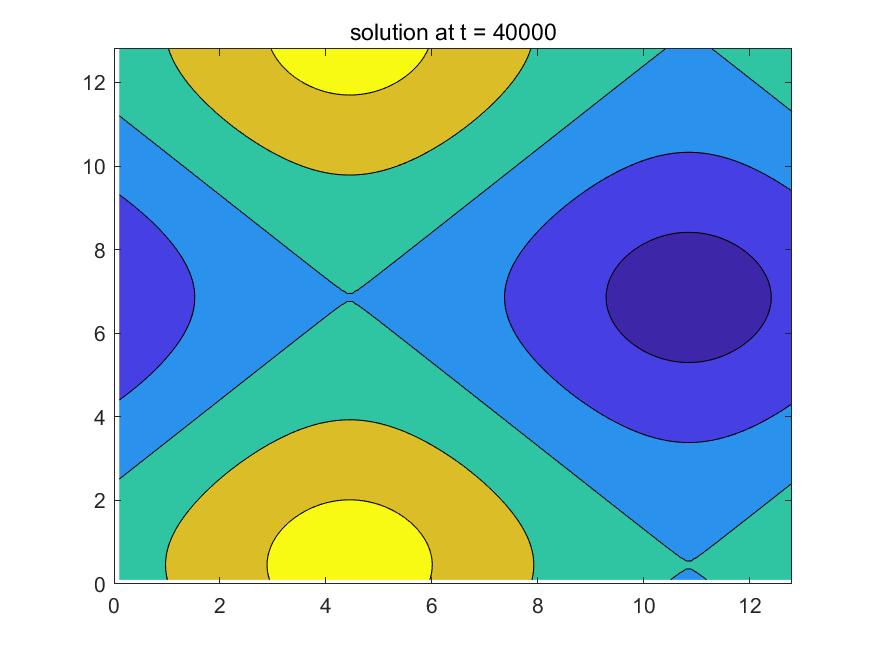}
		\end{minipage}
	 }
	\caption{Snapshots of the numerical solutions for scheme \eqref{eq:second order nss}.}\label{fig: ss}
\end{figure}

\begin{figure}[ht]
	\begin{center}
		 \includegraphics[width =0.8\textwidth]{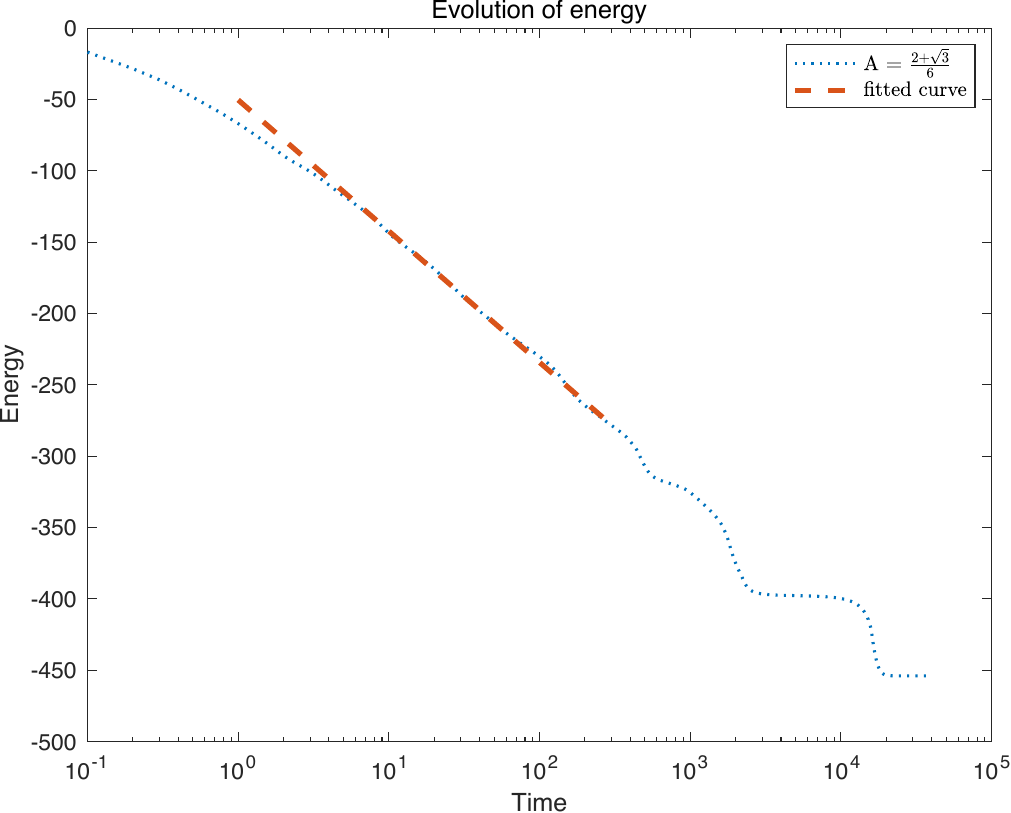}
		\caption{Semi-log plot of the energy $E$. Fitted line has the form $a \ln (t) + b$, with coefficients $a = -39.93$, $b = -50.33$.} \label{fig: energy}
	\end{center}
\end{figure}

\begin{figure}[ht]
	\begin{center}
		 \includegraphics[width=0.4\textwidth]{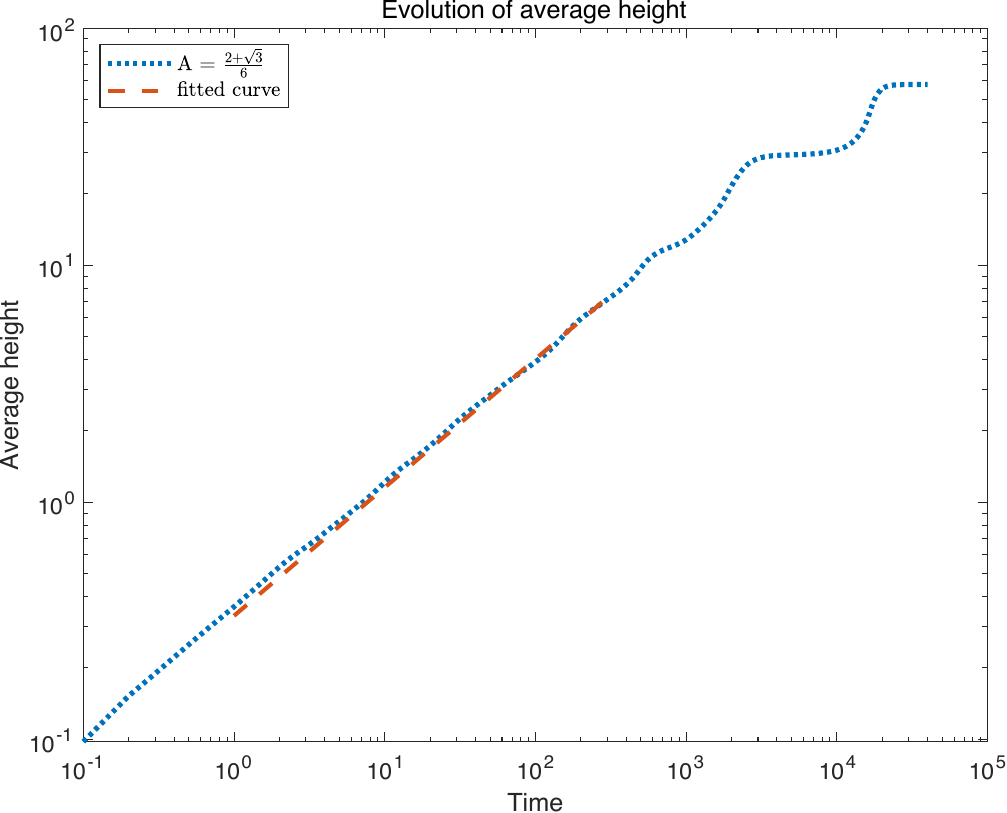}\;
		 \includegraphics[width=0.4\textwidth]{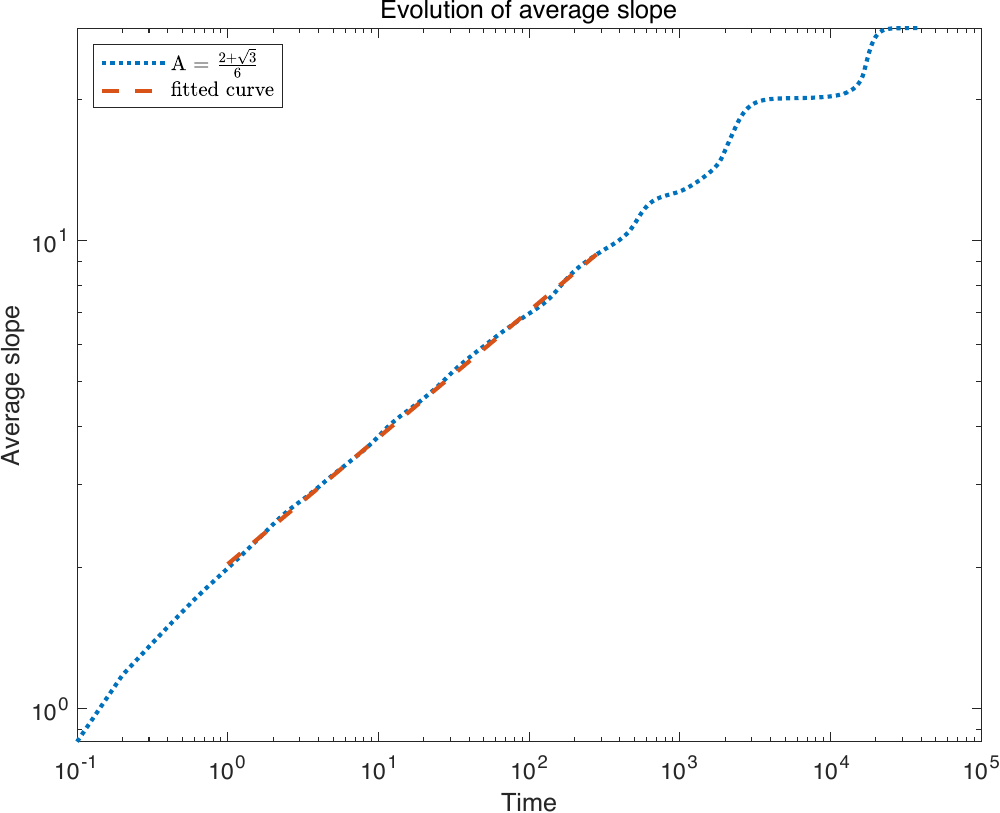}
		\caption{The log-log plot of (1) the average surface height $h$ and (2) the average slope $m$. Fitted lines have the form $a t^b$, with coefficients (1) $a = 0.3319$, $b = 0.5414$ and (2) $a = 2.032$, $b = 0.2705$.} \label{fig: roughness and slope}
	\end{center}
\end{figure}

In this subsection, we simulate the physically interesting coarsening process. The  parameters are $L=4\pi$, $\epsilon=0.005$, $T=40000$. The variable-step time nodes are made upon a uniform grid $\Delta t=0.001$ plus 10\% random perturbation. Some relevant physical quantities, the energy $E$, the average height $h$ and the average slope $m$ are defined by
\begin{align*} 
E(u) &=\left(-\frac{1}{2} \ln \left(1+|\nabla u|^{2}\right), 1\right)+\frac{\varepsilon}{2}\|\Delta u\|^{2}, \\ 
h(u, t) &=\sqrt{\frac{h^{2}}{|\Omega|} \sum_{x}|u(x, t)-\bar{u}(t)|^{2}}, \quad \text { with } \quad \bar{u}(t):=\frac{h^{2}}{|\Omega|} \sum_{x} u(x, t), \\
m(u, t) &=\sqrt{\frac{h^{2}}{|\Omega|} \sum_{x}\left|\nabla u\left(\mathbf{x}_{i, j}, t\right)\right|^{2}},
\end{align*}
and it has been proved $E\sim O(-\ln(t)), \ h\sim O(t^\frac12)$ and $m\sim O(t^\frac14)$ as $t\rightarrow \infty$ (see \cite{golubovic1997interfacial, li2003thin, li2004epitaxial}). The snapshots of numerical solution at time $t=0.99997,\ 5000, \ 15000,\ 20000,\ 30000, \ 40000$ are shown in Figure~\ref{fig: ss}. The scaling laws of energy $E$, average height $h$ and average slope $m$ are verified in Figures~\ref{fig: energy}--\ref{fig: roughness and slope}.

\subsection{Adaptive results}

The solutions of phase field models can exhibit rapid variations over short time intervals, while remaining relatively steady over others. One key benefit of unconditionally energy stable schemes is their natural compatibility with adaptive time-stepping algorithms, where the time step is determined solely by accuracy requirements rather than stability constraints. In contrast, many other numerical schemes face significant challenges when combined with adaptive time stepping, as they often lack robust stability under variable time steps. This highlights the importance of high-order unconditionally energy stable variable step methods.

\begin{figure}[ht]
    \centering
     \includegraphics[width=1\linewidth]{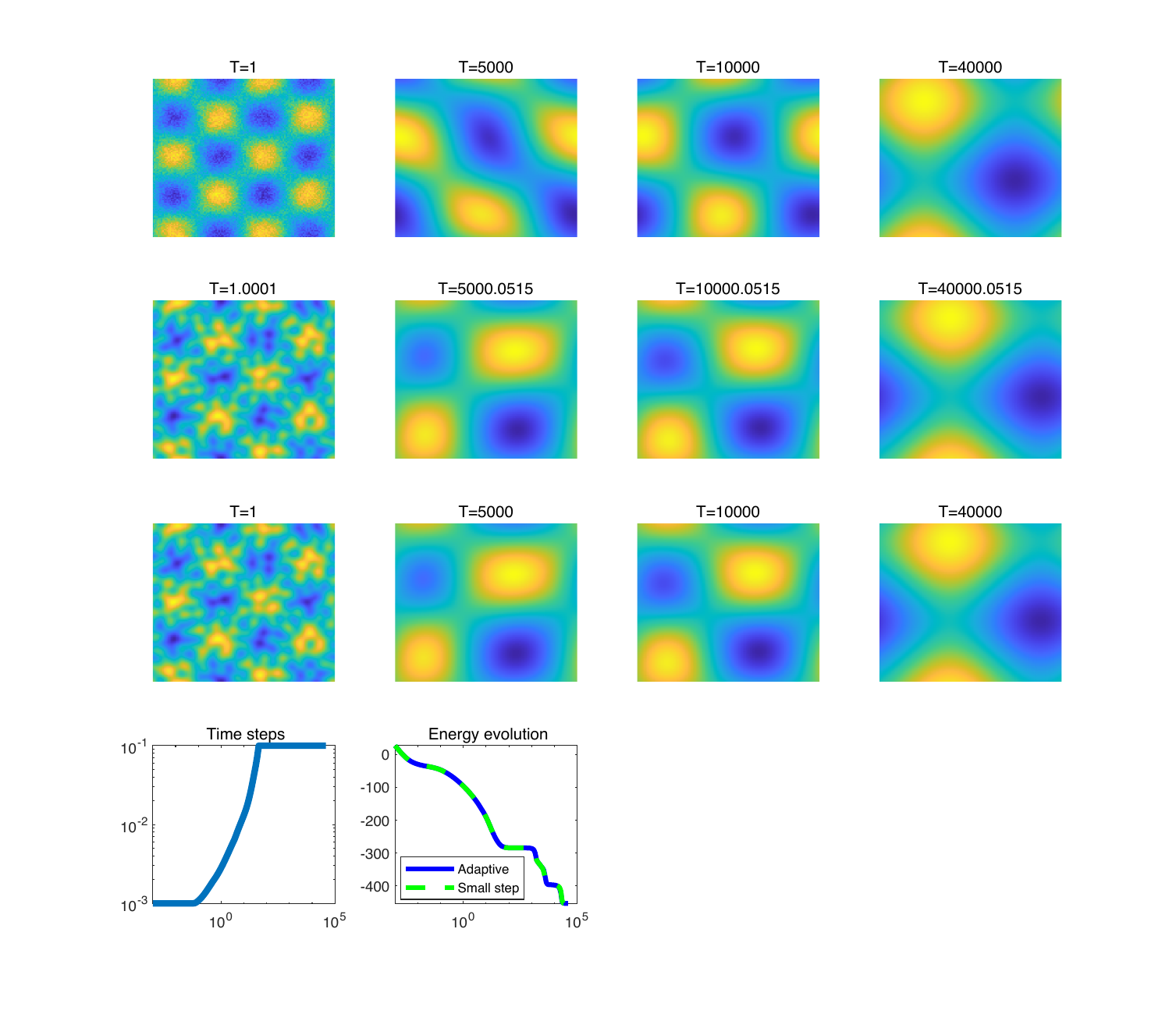}
    \caption{Solutions using large time steps $\tau=0.1$ (first line), adaptive time steps (second line) and small time steps $\tau=0.001$ (third line), and time steps curve and energy evolution (fourth line)}
    \label{fig:adaptivefig}
\end{figure}

For gradient flows, there are several effective adaptive time stepping strategies, see \cite{BDF2CHvarstep,largestepCH,adapCH}. Here we use the strategy summarized in the following Algorithm, where the time step is updated by the formula 
\begin{equation}
    A_{dp} (e,\tau)=\rho \left(\frac{tol}{e}\right)^{r} \tau,
\end{equation}
along with the restriction of the minimum and maximum time steps. In the above formula, $\rho$ is a default safety coefficient, $tol$ is a reference tolerance, $e$ is the relative error computed at each time level in Step 3, and $r$ is the adaptive rate. In our numerical examples, we set $\rho=0.95$ and $tol=10^{-3}$, the minimum time step is $10^{-3}$, while the maximum time step is $\tau=10^{-1}$ for Figure \ref{fig:adaptivefig}. The initial condition is given by $u=\sin(x)\cos(y)+0.5*(2*\text{rand}-1)$ where $\text{rand}$ represents a uniformly distributed random noise in $[0,1]$, and the initial time step is taken as the minimum time step.

\begin{algorithm}[ht]
\caption{Adaptive time stepping procedure}
\textbf{Given:} $U^n,\tau_n$\\
\textbf{Step 1.} Compute $U^{n+1}_1$ by the first-order ETD scheme with $\tau_n$.\\
\textbf{Step 2.} Compute $U^{n+1}_2$ by the second-order ETDMS scheme (6.1) with $\tau_n$.\\
\textbf{Step 3.} Calculate $e_{n+1}=\frac{\|U^{n+1}_1-U^{n+1}_2\|}{\|U^{n+1}_2\|}$.\\
\textbf{Step 4.} If $e_{n+1}>tol$, recalculate the time step $\tau_n\xleftarrow[]{} \max\{\tau_{min},\min\{A_{dp}(e_{n+1},\tau_n),\tau_{max}\}\}$,\\
\textbf{Step 5.} goto Step 1.\\
\textbf{Step 6.} else, update the time step $\tau_{n+1}\xleftarrow[]{}\max\{\tau_{min},\min\{A_{dp}(e_{n+1},\tau_n),\tau_{max}\}\}$.\\
\textbf{Step 7.} endif
\end{algorithm}

We take $\epsilon=0.005$, $N=128$ and as a comparison, the first and third rows in Fugure \ref{fig:adaptivefig} exhibit snapshots of the numerical result generated by the second-order ETDMS solutions with uniform time step size  $\tau=10^{-1}$ and $\tau=10^{-3}$, respectively, while the second row displays snapshots generated by the adoptive scheme starting from the same initial data.  The last row presents the magnitudes of time steps and the energy evolution. We observe that the solution using large time step (first row) even cannot get the correct topological changes. We also observe from the last row that the adaptive time steps basically stay around the maximum size at large time, indicating that the computational cost is almost the same as the large-step solution. In addition, we can also observe that the energy curve has a similar evolution as in Figure 2. This experiment proves that with the unconditionally energy-stable ETDMS schemes, the adaptive algorithm only takes as little computational cost as that of large time steps, while achieving the same level of accuracy as small time steps.

\section{Conclusion} \label{sec:conclusion}
Extending our recent work \cite{chen2021ETDMS}, a generic variable-step $k^{th}$-order linear scheme was proposed for gradient flows by using ETD-MS method. Under the assumption that the nonlinearity is Lipschitz continuous in some appropriate sense, we have demonstrated that the scheme is long time stable in the sense that a modified energy is a monotonic function of the time. In addition, we proved the $k^{th}$-order accuracy in the $\ell^{\infty}(0,T; V^\lambda)$ norm  assuming a maximal time step-size $1$ {and a bound on local (neighboring-$k$) time-step ratio}. Numerical experiments on the thin film epitaxial growth without slope selection model confirm the stability and optimal rate of convergence. {Using local error as an adaptive strategy, we also demonstrate the potential of the scheme as an efficient adaptive-in-time algorithm.}

\section*{Acknowledgements} 
	This work is supported in part by the grants NSFC 12241101, NSFC 12471369 (W. Chen), NSFC 11871159, NSFC 12271237 (X.Wang). 

\bibliographystyle{99}

% \begin{thebibliography}{99}

% \bibitem{1} B. D. O. Anderson and J. B. Moore, {\it Linear optimal
% control},  Prentice-Hall,  Englewood Cliffs, NJ, 1971.

% \bibitem{2}S. Aubry and P.Y. Le Daeron, {\it The discrete
% Frenkel-Kontorova model and its extensions I},  Physica D  {\bf 8} (1983), 381--422.

% \bibitem{3} J. Baumeister, A. Leitao and G. N. Silva, {\it On
% the value function for nonautonomous optimal control problem with
% infinite horizon},  Systems Control Lett. {\bf 56}  (2007),  188--196.

% \bibitem{4} J. Blot, {\it Infinite-horizon Pontryagin
% principles without invertibility},  J. Nonlinear Convex Anal.
% {\bf 10} (2009),  177--189.

% \bibitem{5} J. Blot and P. Cartigny, {\it  Optimality in
% infinite-horizon variational problems under sign conditions},
% J. Optim. Theory Appl. {\bf 106} (2000),  411--419.

% \end{thebibliography}

\end{document}